\documentclass[12pt]{article}

\usepackage{latexsym,amsmath}
\usepackage{times}
\usepackage{amssymb}

\newtheorem{thm}{Theorem}[section]

\newtheorem{lemma}[thm]{Lemma}

\newtheorem{conj}[thm]{Conjecture}

\newcommand{\beq}[1]{\begin{equation}\label{#1}}
\newcommand{\enq}[0]{\end{equation}}

\newcommand{\qed}[0]{{\hspace*{\fill}\mbox{$\Box$}}}

\newcommand{\cI}[0]{{\cal I}}

\newcommand{\ga}[0]{\alpha}

\newcommand{\gl}[0]{\lambda}

\begin{document}
\renewcommand{\thefootnote}{\fnsymbol{footnote}}

\title{An upper bound for the number of independent sets in regular graphs}

\author{David Galvin}

\author{David Galvin\thanks{Department of Mathematics,
University of Notre Dame, 255 Hurley Hall, Notre Dame IN
46556; dgalvin1@nd.edu; + 1 574 631 7245 (telephone); + 1 574 631 6579 (fax)}}

\maketitle

\begin{abstract}
Write ${\cal I}(G)$ for the set of independent sets of a graph $G$ and $i(G)$ for $|{\cal I}(G)|$. It has been conjectured (by Alon and Kahn) that for an $N$-vertex, $d$-regular graph $G$,
$$
i(G) \leq \left(2^{d+1}-1\right)^{N/2d}.
$$
If true, this bound would be tight, being achieved by the disjoint union of $N/2d$ copies of $K_{d,d}$. Kahn established the bound for bipartite $G$, and later gave an argument that established
$$
i(G)\leq 2^{\frac{N}{2}\left(1+\frac{2}{d}\right)}
$$
for $G$ not necessarily bipartite. In this note, we improve this to
$$
i(G)\leq 2^{\frac{N}{2}\left(1+\frac{1+o(1)}{d}\right)}
$$
where $o(1) \rightarrow 0$ as $d \rightarrow \infty$, which matches the conjectured upper bound in the first two terms of the exponent.

We obtain this bound as a corollary of a new upper bound on the independent set polynomial $P(\gl,G)=\sum_{I \in {\cal I}(G)} \gl^{|I|}$ of an $N$-vertex, $d$-regular graph $G$, namely
$$
P(\gl,G) \leq (1+\gl)^{\frac{N}{2}} 2^{\frac{N(1+o(1))}{2d}}
$$
valid for all $\gl > 0$. This also allows us to improve the bounds obtained recently by Carroll, Galvin and Tetali on the number of independent sets of a fixed size in a regular graph.

\end{abstract}

\section{Introduction}

For a (simple, finite, undirected) graph $G$ write ${\cal I}(G)$ for the set of independent sets of $G$ (sets of vertices no two of which are adjacent) and $i(G)$ for $|{\cal I}(G)|$. How large can $i(G)$ be? For the class of $N$-vertex, $d$-regular graphs, this question has received some attention, with a succession of bounds having appeared going back to the early 1990's.

A trivial upper bound is $i(G) \leq 2^N$. In \cite{Alekseev}, Alekseev gave the first non-trivial bound, establishing $i(G) \leq 3^\frac{N}{2}$.
In fact he showed that for any graph $G$ on $N$ vertices (not necessarily regular)
\begin{equation} \label{Alekseev-bound}
i(G) \leq \left(1+\frac{N}{\alpha}\right)^\alpha
\end{equation}
where $\alpha=\alpha(G)$ is the size of the largest independent set in $G$; note that for $N$-vertex, $d$-regular $G$, $\alpha(G)\leq \frac{N}{2}$. Around the same time, in the process of resolving a question of Erd\H{o}s and Cameron on sum-free sets in Abelian groups, Alon \cite{Alon} substantially improved this to
\begin{equation} \label{Alon-bound}
i(G) \leq \exp_2\left\{\frac{N}{2}\left(1+ \frac{C}{d^{1/10}}\right)\right\}
\end{equation}
for some constant $C>0$.
This bound is best possible in the leading term of the exponent: the graph $\frac{N}{2d}K_{d,d}$ consisting of a disjoint union of $N/2d$ copies of $K_{d,d}$ satisfies
$$
i\left(\frac{N}{2d}K_{d,d}\right) = \left(2^{d+1}-1\right)^\frac{N}{2d} =  \exp_2\left\{\frac{N}{2}\left(1+ \frac{1}{d} - \frac{1 + o(1)}{(2\ln 2)d2^d}\right)\right\}
$$
where $o(1) \rightarrow 0$ as $d \rightarrow \infty$.
Alon speculated that perhaps among all $N$-vertex, $d$-regular graphs, $\frac{N}{2d}K_{d,d}$ is the one that admits the greatest number of independent sets.
\begin{conj} \label{conj-Alon-Kahn}
For any $N$-vertex, $d$-regular graph $G$,
$$
i(G) \leq \left(2^{d+1}-1\right)^\frac{N}{2d}.
$$
\end{conj}
Kahn \cite{Kahn} used entropy methods to prove the above bound for $N$-vertex, $d$-regular bipartite graphs, and in the same paper formally conjectured that the bound should hold for all graphs.

Further progress was made by Sapozhenko \cite{Sap5}, who used a very simple counting argument to improve (\ref{Alon-bound}) to
$$
i(G) \leq \exp_2\left\{\frac{N}{2}\left(1+ C\sqrt{\frac{\log d}{d}}\right)\right\}
$$
for some constant $C>0$. (In this note ``$\log$'' will always indicate the base $2$ logarithm).

The next substantial improvement was made by Kahn (personal communication to the author; the proof appears in \cite{MadimanTetali}, where it is generalized to the context of graph homomorphisms), who obtained
\begin{equation} \label{inq-kahn-nonbipartite}
i(G)\leq \exp_2\left\{\frac{N}{2}\left(1+\frac{2}{d}\right)\right\}.
\end{equation}

The aim of this note is to improve (\ref{inq-kahn-nonbipartite}) to the following.
\begin{thm} \label{thm-galvin-nonbipartite}
There is a constant $C>0$ such that for any $d$-regular, $N$-vertex graph $G$,
$$
i(G) \leq \exp_2\left\{\frac{N}{2}\left(1+\frac{1}{d}+\frac{C}{d}\sqrt{\frac{\log d}{d}}\right)\right\}.
$$
\end{thm}
This matches the first two terms in the exponent of $i\left(\frac{N}{2d}K_{d,d}\right)$.
The proof of Theorem \ref{thm-galvin-nonbipartite} combines the idea used to prove (\ref{inq-kahn-nonbipartite}) with a recent theorem of Sapozhenko bounding the number of independent sets in a regular graph in terms of the size of the largest independent set (see Lemma \ref{thm-sap-maxalpha-weighted}). The basic idea is to treat two cases. If $G$ has no large independent sets, then Sapozhenko's result implies that it has
few independent sets.
On the
other hand if $G$ has a large independent set then it is close to bipartite (the case for which Conjecture \ref{conj-Alon-Kahn} is resolved) and the method used to prove (\ref{inq-kahn-nonbipartite}) can be modified to exploit this fact and obtain a bound closer to that of Conjecture \ref{conj-Alon-Kahn}.

\medskip

What we actually prove is a weighted generalization of Theorem \ref{thm-galvin-nonbipartite}. The independent set (or stable set) polynomial of $G$ (first introduced explicitly by Gutman and Harary \cite{GutmanHarary}) is defined as
$$
P(\gl,G) = \sum_{I \in {\cal I}(G)} \gl^{|I|}.
$$
This is also referred to as the partition function of the independent set (or hard-core) model on $G$ with activity $\gl$. In \cite{GalvinTetali} the analog of Conjecture \ref{conj-Alon-Kahn} was obtained for $N$-vertex, $d$-regular bipartite $G$:
$$
P(\gl,G) \leq (2(1+\gl)^d -1)^{\frac{N}{2d}}~\left(=P\left(\gl, K_{d,d}\right)^{\frac{N}{2d}}\right)
$$
for all $\gl >0$ (the case $\gl \geq 1$ was already dealt with in \cite{Kahn}), and it was conjectured that this bound should hold for non-bipartite $G$ also. In \cite{CarrollGalvinTetali} the analog of (\ref{inq-kahn-nonbipartite}) was obtained for $N$-vertex, $d$-regular $G$:
\begin{equation} \label{g-weighted-bound}
P(\gl, G) \leq (1+\gl)^{\frac{N}{2}} 2^{\frac{N}{d}}.
\end{equation}
By employing a weighted generalization of (\ref{Alekseev-bound}) (see Lemma \ref{lem-Alekseev-weighted}) we improve (\ref{g-weighted-bound}) to the following.
\begin{thm} \label{thm-galvin-nonbipartite-weighted}
For all $\gl > 0$ there is a constant $C_\gl>0$ such that for $d$-regular, $N$-vertex $G$,
$$
P(\gl, G) \leq (1+\gl)^{\frac{N}{2}}\exp_2\left\{\frac{N}{2d}\left(1+C_\gl\sqrt{\frac{\log d}{d}}\right)\right\}.
$$
\end{thm}
Note that this reduces to Theorem \ref{thm-galvin-nonbipartite} when $\gl=1$.

\medskip

Theorem \ref{thm-galvin-nonbipartite-weighted} has consequences for the number of independent sets of a fixed size in a regular graph. With regards to this, Kahn \cite{Kahn} made the following conjecture. Here $i_t(G)$ is the number of independent sets in $G$ of size $t$.
\begin{conj} \label{conj-fixed-size}
If $G$ is an $N$-vertex, $d$-regular graph with $2d|N$, then for each $0 \leq t \leq N$,
$$
i_t(G) \leq i_t\left(\frac{N}{2d}K_{d,d}\right).
$$
\end{conj}
In \cite{CarrollGalvinTetali}, asymptotic evidence is provided for this conjecture in the sense that if $N$, $d$ and $t$ are sequences satisfying $t=\alpha
N/2$ for some fixed $\alpha \in (0,1)$ and $G$ is a sequence
of $N$-vertex, $d$-regular graphs, then
\begin{equation} \label{GCT}
i_t(G) \leq \left\{
\begin{array}{ll}
\exp_2\left\{\frac{N}{2}\left(H\left(\alpha\right)+\frac{2}{d}\right)\right\} &
\mbox{in general and} \\
& \\
\exp_2\left\{\frac{N}{2}\left(H\left(\alpha\right)+\frac{1}{d}\right)\right\} &
\mbox{if $G$ is bipartite,}
\end{array}
\right.
\end{equation}
where $H(\cdot)$ is the binary entropy function. On the other hand,
if $N=\omega(d \log d)$ and $d=\omega(1)$ then
$$
i_t\left(\frac{N}{2d}K_{d,d}\right)  \geq
\exp_2\left\{\frac{N}{2}\left(H\left(\alpha\right)+\frac{1-o(1)}{d}\right)\right\}
$$
(all as $d \rightarrow \infty$). Using Theorem \ref{thm-galvin-nonbipartite-weighted} in place of (\ref{g-weighted-bound}) in the derivation of the first bound in (\ref{GCT}) we get the immediate improvement that there is a constant $c_\ga>0$ such that for all $N$-vertex, $d$-regular $G$
$$
i_t(G) \leq \exp_2\left\{\frac{N}{2}\left(H\left(\alpha\right)+\frac{1}{d} + \frac{c_\ga}{d}\sqrt{\frac{\log d}{d}}\right)\right\},
$$
so that the upper bound for general $G$ matches the conjectured bound in the first two terms of the exponent (in the range $N=\omega(d \log d)$, $d=\omega(1)$ and $t=\alpha N/2$).

\medskip

Section \ref{sec-tools} gives the three lemmas that we need for the proof of Theorem \ref{thm-galvin-nonbipartite-weighted}, while the proof of the theorem is given in Section \ref{sec-proofs-1}.

\section{Tools} \label{sec-tools}

We begin by recalling a result from \cite{CarrollGalvinTetali} which is a slight refinement of (\ref{g-weighted-bound}) (see the derivation of (10) in that reference). For a total order $\prec$ on $V(G)$ and for each $v \in V(G)$ write $P_\prec(v)$ for $\{w \in V(G):\{w,v\} \in E(G),~ w \prec v\}$
and $p_\prec(v)$ for $|P_\prec(v)|$. Note that $\sum_{v \in V(G)} p_\prec(v) = |E(G)|$ ($=Nd/2$ when $G$ is $N$-vertex and $d$-regular).
\begin{lemma} \label{thm-madtet-weighted}
For any $d$-regular, $N$-vertex graph $G$ and any total order $\prec$ on $V(G)$,
$$
P(\gl, G) \leq \prod_{v \in V(G)}
\left(2(1+\gl)^{p_\prec(v)}-1\right)^\frac{1}{d}.
$$
\end{lemma}

\medskip

Next, we give a weighted generalization of (\ref{Alekseev-bound}).
\begin{lemma} \label{lem-Alekseev-weighted}
For any $N$-vertex graph $G$ (not necessarily regular) with $\alpha(G)=\alpha$, and any $\gl > 0$,
$$
P(\gl, G) \leq \left(1+\frac{\gl N}{\alpha}\right)^\alpha
$$
with equality if and only if $G$ is the disjoint union of complete graphs all of the same order.
\end{lemma}

\medskip

\noindent {\em Proof}: We follow closely the proof of (\ref{Alekseev-bound}) that appears in \cite{Alekseev2}, making along the way the changes needed to introduce $\gl$.

We first observe that it is enough to prove the bound for connected $G$. Indeed, if $G$ has components $C_1, \ldots, C_n$ with $|C_i|=N_i$ and $\alpha(G[C_i])=\alpha_i$ then  $\sum N_i = N$ and $\sum \alpha_i = \alpha$. Using Jensen's inequality for (\ref{Jensen}) we have
\begin{eqnarray}
P(\gl, G) & = & \prod_{i=1}^n P(\gl, G[C_i]) \nonumber \\
& \leq & \prod_{i=1}^n \left(1+\frac{\gl N_i}{\alpha_i}\right)^{\alpha_i} \nonumber \\
& \leq & \left( \frac{\sum_{i=1}^n \alpha_i\left(1+\frac{\gl N_i}{\alpha_i}\right)}{\alpha} \right)^\alpha \label{Jensen} \\
& = & \left(1+\frac{\gl N}{\alpha}\right)^\alpha, \nonumber
\end{eqnarray}
with equality in (\ref{Jensen}) if and only if all the $N_i$ are equal (i.e., all components of $G$ have the same order).

We prove the bound for connected $G$ by induction on the number of vertices, with the single-vertex case trivial. For connected $G$ with more than one vertex, let $v$ be a vertex of maximum degree $\Delta$. We use the recurrence
$$
P(\gl, G) = P(\gl, G-v) + \gl P(\gl, G-v-N(v))
$$
(where $N(v)$ is the set of vertices in $G$ adjacent to $v$), which follows from the fact that there is a bijection from independent sets of size $t$ in $G$ which do not contain $v$ to independent sets of size $t$ in $G-v$, and a bijection from independent sets of size $t$ in $G$ which do contain $v$ to independent sets  of size $t-1$ in $G-v-N(v)$. Since $\alpha(G-v) \leq \alpha$ and $\alpha(G-v-N(v)) \leq \alpha-1$ we have by induction and the fact that $\left(1+\frac{a}{x}\right)^x$ is increasing in $x>0$ for all $a>0$
\begin{equation} \label{inq-using-recurrence}
P(\gl, G) \leq \left(1+\frac{\gl(N-1)}{\alpha}\right)^\alpha + \gl \left(1+\frac{\gl(N-\Delta-1)}{\alpha-1}\right)^{\alpha - 1}.
\end{equation}
We upper bound the right-hand side of (\ref{inq-using-recurrence}) by observing that $\Delta \geq \frac{N-1}{\alpha}$. For $G$ complete or an odd cycle, this is immediate, and for all other connected $G$ the stronger bound $\Delta \geq \frac{N}{\alpha}$ follows from Brooks' theorem.
Inserting into (\ref{inq-using-recurrence}) we obtain
$$
P(\gl, G)  \leq   \left(1+\frac{\gl(N-1)}{\alpha}\right)^\alpha + \gl \left(1+\frac{\gl(N-1)}{\alpha}\right)^{\alpha - 1}
$$
and so
\begin{eqnarray}
\frac{P(\gl, G) - \left(1+\frac{\gl N}{\alpha}\right)^\alpha}{\left(1+\frac{\gl(N-1)}{\alpha}\right)^\alpha} & \leq &
1 + \frac{\gl \alpha}{\alpha + \gl(N-1)} - \left(1+\frac{\gl}{\alpha + \gl(N-1)}\right)^\alpha \nonumber \\
& \leq & 0 \nonumber
\end{eqnarray}
with equality if and only if $\alpha=1$ (i.e., $G$ is a complete graph). \qed

\medskip

Finally, we give a weighted variant of a recent result of Sapozhenko \cite{Sap4} bounding the number of independent sets in a regular graph in terms of $\alpha(G)$.
\begin{lemma} \label{thm-sap-maxalpha-weighted}
There is a constant $c>0$ such that for any $d$-regular, $N$-vertex graph $G$ with $d \geq 2$ and $\alpha(G)=\alpha$ and any $\gl>0$,
$$
P(\gl,G) \leq  \left(1+\frac{\gl N}{2\alpha}\right)^\alpha \exp_2\left\{c N\sqrt{\frac{\log d}{d}}\right\}.
$$
\end{lemma}

\medskip

\noindent {\em Proof}: We follow the proof from \cite{Sap4} of the case $\gl = 1$, replacing an appeal to (\ref{Alekseev-bound}) in that proof with an appeal to Lemma \ref{lem-Alekseev-weighted}.

Fix an integer $0 < \varphi < d$. For an independent set $I \in \cI(G)$, recursively construct
sets $T(I)$ and $D(T)$ as follows. Pick $u_1 \in I$ and set $T_1 =
\{u_1\}$. Given $T_m =\{u_1, \ldots, u_m\}$, if there is $u_{m+1}
\in I$ with $N(u_{m+1})\setminus N(T_m) \geq \varphi$, then set $T_{m+1}=\{u_1, \ldots,
u_{m+1}\}$. If there is no such $u_{m+1}$, then set $T=T_m$ and
$$
D(T)=\{v \in V(G)\setminus N(T): N(v)\setminus N(T) < \varphi\}.
$$
Note that
\begin{equation} \label{Tbound}
|T| \leq \frac{N}{\varphi},
\end{equation}
since each step in the construction of
$T$ removes at least $\varphi$ vertices from consideration; that
\begin{equation} \label{Dfact}
I \subseteq D
\end{equation}
since if $I \setminus D \neq \emptyset$, the
construction of $T$ would not have stopped (note that $N(T) \cap I =
\emptyset$); and that
\begin{equation} \label{Dbound}
|D| \leq \frac{Nd}{2d-\varphi}.
\end{equation}
To see (\ref{Dbound}), consider the bipartite graph with partition classes $D$
and $N(T)$ and edges induced from $G$. This graph has at most $d|N(T)| \leq d(N-|D|)$ edges
(since each vertex in $N(T)$ has at most $d$ edges to $D$, and there
are at most $N-|D|$ such vertices), and at least $(d-\varphi)|D|$ edges
(since each vertex in $D$ has at least $d-\varphi$ edges to $N(T)$).
Putting these two inequalities together gives (\ref{Dbound}).

Combining (\ref{Tbound}), (\ref{Dfact}) and (\ref{Dbound}) we see that we can construct all $I \in \cI(G)$ by first picking a $T \subseteq V$ of size at most $N/\varphi$, next constructing $D(T)$, and finally generating all independent sets of the subgraph of $G$ induced by $D(T)$. This graph inherits from $G$ the property that all  independent sets have size at most $\alpha$, so using Lemma \ref{lem-Alekseev-weighted} it follows that
\begin{equation} \label{int3}
P(\gl, G) \leq \sum_{t \leq N/\varphi} {N \choose t} \left(1+\frac{\gl Nd}{(2d-\varphi)\alpha}\right)^\alpha.
\end{equation}
We bound
$$
1+\frac{\gl Nd}{(2d-\varphi)\alpha} = 1 + \frac{\gl N}{2\alpha} \left(\frac{2d}{2d-\varphi}\right) \leq \left(1 + \frac{\gl N}{2\alpha}\right) \left(\frac{2d}{2d-\varphi}\right)
$$
so that
\begin{equation} \label{int4}
\left(1+\frac{\gl Nd}{(2d-\varphi)\alpha}\right)^\alpha \leq \left(1+\frac{\gl N}{2\alpha}\right)^\alpha \left(\frac{2d}{2d-\varphi}\right)^N
\end{equation}
and naively bound
\begin{equation} \label{int5}
\sum_{t \leq N/\varphi} {N \choose t} \leq \left(\frac{N}{\varphi}+1\right) {N \choose N/\varphi} \leq \exp_2\left\{\frac{3N\log (e\varphi)}{\varphi}\right\}
\end{equation}
(the second inequality mainly using ${n \choose k} \leq (en/k)^k$).
Taking $\varphi = [\sqrt{d \log d}]$ the bounds in (\ref{int3}), (\ref{int4}) and (\ref{int5}) combine to yield
$$
P(\gl, G) \leq  \left(1+\frac{\gl N}{2\alpha}\right)^\alpha \exp_2\left\{c N\sqrt{\frac{\log d}{d}}\right\}
$$
for some $c>0$, as claimed.
\qed

\section{Proof of Theorem \ref{thm-galvin-nonbipartite-weighted}} \label{sec-proofs-1}

We assume throughout that $d \geq 2$, since the theorem is straightforward for $d=1$. We begin by considering those $G$ for which
$$
\alpha(G) \geq \frac{N}{2}\left(1-C_\gl\sqrt{\frac{\log d}{d}}\right),
$$
where $C_\gl$ is a constant that will be determined later. In this case we use Lemma \ref{thm-madtet-weighted}.
For each $v$ with $p_\prec(v)=0$ we have $2(1+\gl)^{p_\prec(v)}-1=1$ and so if $\prec$ begins by listing the vertices of an independent set $I$ then
\begin{eqnarray}
P(\gl, G) & \leq & \prod_{v \in V(G)\setminus I} \left(2(1+\gl)^{p_\prec(v)}-1\right)^\frac{1}{d} \nonumber \\
& \leq & (1+\gl)^{\frac{1}{d}\sum_{v \in V(G)\setminus I} p_\prec(v)}\exp_2\left\{\frac{|V(G)\setminus I|}{d}\right\} \nonumber \\
& = & (1+\gl)^\frac{N}{2}\exp_2\left\{\frac{N-|I|}{d}\right\} \label{int2}.
\end{eqnarray}
Choosing $I$ to be an independent set of size $\alpha(G)$ we get from (\ref{int2}) that
\begin{equation} \label{case-largealpha-weighted}
P(\gl, G) \leq (1+\gl)^{\frac{N}{2}} \exp_2\left\{\frac{N}{2d}\left(1+C_\gl\sqrt{\frac{\log d}{d}}\right)\right\}.
\end{equation}

\medskip

We use Lemma \ref{thm-sap-maxalpha-weighted} to bound $P(\gl, G)$ in the case when $\alpha(G)$ satisfies
$$
\alpha(G) < \frac{N}{2}\left(1-C_\gl\sqrt{\frac{\log d}{d}}\right).
$$
For typographical convenience, in what follows we write $x$ for $C_\gl\sqrt{\log d/d}$. Since $\left(1+\frac{\gl N}{2\alpha}\right)^\alpha$ is increasing in $\alpha > 0$ for $\gl, N > 0$ we have
\begin{eqnarray*}
\left(1+\frac{\gl N}{2\alpha}\right)^\alpha & \leq & \left(1+\frac{\gl}{1-x}\right)^{\frac{N(1-x)}{2}} \\
& = & (1+\gl)^\frac{N}{2} \left(1+\frac{x\gl}{(1-x)(1+\gl)}\right)^{\frac{N(1-x)}{2}} (1+\gl)^{-\frac{Nx}{2}} \\
& \leq & (1+\gl)^\frac{N}{2} \exp\left\{\frac{Nx}{2}\left(\frac{\gl}{1+\gl}-\ln(1+\gl)\right)\right\}.
\end{eqnarray*}
Since $\ln(1+\gl)>\frac{\gl}{1+\gl}$ for all $\gl > 0$ the exponent above is negative. By choosing $C_\gl> 0$ to satisfy
$$
\left(\ln(1+\gl) - \frac{\gl}{1+\gl}\right)\frac{C_\gl}{2\ln 2} = c
$$
where $c$ is the constant appearing in the bound in Lemma \ref{thm-sap-maxalpha-weighted}, we have
\begin{equation} \label{case-smallalpha-weighted}
P(\gl, G) \leq (1+\gl)^\frac{N}{2}
\end{equation}
in this case. Combining (\ref{case-largealpha-weighted}) and (\ref{case-smallalpha-weighted}) we obtain
$$
P(\gl, G) \leq (1+\gl)^{\frac{N}{2}} \exp_2\left\{\frac{N}{2d}\left(1+C_\gl\sqrt{\frac{\log d}{d}}\right)\right\}
$$
for all $G$, completing the proof of Theorem \ref{thm-galvin-nonbipartite-weighted}.

\end{document}